\documentclass[12pt]{article}

\usepackage{epsf,graphicx,latexsym}

\newcommand{\dst} {\displaystyle}

\newcommand {\Hol}{\mathop{\rm Hol}\nolimits}

\renewcommand {\Re}{\mathop{\rm Re}\nolimits}

\newcommand {\DD}{\Delta}

\newcommand {\gmb}{\mathop{\rm G}\nolimits(\mu,\beta)}
\newcommand {\gub}{\mathop{\rm G}\nolimits(1,\beta)}

\newfont{\bbb}{msbm10 at 12pt}
\def\Bbb#1{\hbox{\bbb #1}}

\newcommand{\C}{{\Bbb C}}
\newcommand{\R}{{\Bbb R}}

\usepackage{amsfonts}
\usepackage{amssymb}
\usepackage{latexsym}

\newcommand{\pl}{\partial}
\newcommand{\pr}{\noindent{\bf Proof.}\quad }
\newcommand{\epr}{\ $\blacksquare$ \vspace{2mm} }

\newcommand{\be} {\begin{eqnarray}}
\newcommand{\ee} {\end{eqnarray}}
\newcommand{\bep} {\begin{eqnarray*}}
\newcommand{\eep} {\end{eqnarray*}}

\newtheorem{defin}{Definition}[section]
\newtheorem{theorem}{Theorem}[section]
\newtheorem{corol}{Corollary}[section]
\newtheorem{propo}{Proposition}[section]

\begin{document}

\title{Covering and distortion theorems for spirallike functions with respect to a boundary
point\thanks{\footnotesize{\it 2000 Mathematics Subject
Classification}: Primary 30C45}\thanks{\footnotesize{\it Key words
and phrases}: spirallike functions with respect to a boundary
point, distortion theorem, covering theorem}}
\author{Mark Elin \\ {\small Department of Mathematics, ORT  Braude College,} \\
{\small P.O. Box 78, Karmiel 21982, ISRAEL} \\ {\small e-mail:
mark.elin@gmail.com}}
\date{}

\maketitle

\begin{abstract}
In this note we obtain some distortion results for spirallike
functions with respect to a boundary point. In particular, we find
the maximal domain covered by all spirallike functions of order
$\beta$.
\end{abstract}

\setcounter{section}{-1}
\section{Introduction}
\setcounter{equation}{0}

Although the classes of starlike and spirallike functions (taking
zero at zero) were objects of interest during about the last
hundred years, the study of functions having similar geometric
properties with respect to a boundary point was begun only in the
1980's.

A breakthrough in this matter is due to M. S. Robertson
\cite{RMS}, who defined the class of those univalent holomorphic
functions $f\in\Hol(\Delta,\C)$ on the open unit disk $\Delta$
satisfying $f(0)=1$ such that $f(\Delta)$ is starlike with respect
to the boundary point $f(1):=\lim\limits_{r\to1^-}f(r)=0$ and
$f(\Delta)$ lies in a half-plane. Different characterizations of
the class of starlike functions with respect to a boundary point
were obtained in \cite{LA, SH-SEM, E-R-S1, LeA, LA-LA}; properties
of these functions were considered in numerous works (see, for
example, \cite{A-A-S, Z-Q}). Various distortion results for
starlike functions with respect to a boundary point were proved in
\cite{T, C-O, D-L}.

It seems that except for \cite{E-R-S2}, the paper \cite{A-E-S} was
the first where spirallike functions with respect to a boundary
point were studied systematically.

\begin{defin}
A univalent function $f\in\Hol(\Delta,\C)$ normalized by $f(0)=1$
and $f(1):=\lim\limits_{r\to1^-}f(r)=0$ is called spirallike with
respect to a boundary point if there is a number $\mu\in\C,\
\Re\mu>0,$ such that for any point $w\in f(\Delta)$ the curve
$\dst\{e^{-t\mu}w,\ t\ge0\}$ is contained in $f(\Delta)$.

If, in particular, we also have $\mu\in\R$, then $f$ is called
starlike with respect to a boundary point.
\end{defin}

The following assertion was proved in \cite{A-E-S}.

\begin{theorem}\label{AES}
Let $f\in\Hol(\Delta,\C),\ f(0)=1,\ f(1)=0,$ be a spirallike
function with respect to a boundary point. Then there exists a
number $\mu\in\Omega:=\left\{\lambda\in\C:\ |\lambda-1|\le1,\
\lambda\not=0 \right\}$ such that
\be\label{main1}
\Re\left(\frac2\mu\,\frac{zf'(z)}{f(z)}\,+\,\frac{1+z}{1-z}
\right)>0.
\ee
Conversely, if a univalent function $f,\ f(0)=1,\ f(1)=0,$
satisfies (\ref{main1}) for some  $\mu\in\Omega$, then $f$ is a
spirallike function with respect to a boundary point.
\end{theorem}

Further, this class was investigated in \cite{E-S2, SD1, LeA1}. In
particular, combining results of \cite{E-S2} and \cite{A-E-S} one
can formulate the following geometric characterization of images
of spirallike functions with respect to a boundary point.

\begin{theorem}\label{ES}
Let $f\in\Hol(\Delta,\C),\ f(0)=1,\ f(1)=0,$ be a spirallike
function with respect to a boundary point. Then

(i) the following two limits exist finitely:
\begin{equation}\label{mu}
\displaystyle\lim_{r\to1^-}\,\frac{f'(r)(r-1)}{f(r)}\,=\nu
\end{equation}
with $\nu\in\Omega$, and
\begin{equation}\label{a}
\displaystyle\lim_{r\to1^-}\arg f^{1/\nu}(r)=a
\end{equation}
with $\displaystyle|a|<\frac\pi2$;

(ii) inequality (\ref{main1}) holds for $\mu=\nu$ defined by
(\ref{mu}); moreover, any number $\mu$ for which inequality
(\ref{main1}) holds must satisfy $\mu=R\nu,\ R\ge1$;

(iii) the minimal spiral wedge that contains the image $f(\Delta)$
is bounded by two spirals:
\[
w=e^{-\nu t}w_+\mbox{ and }w=e^{-\nu t}w_-,\ t\in\R,
\]
where $\displaystyle w_\pm=\exp\biggl(i\nu(a\pm\frac\pi2)\biggr)$.
\end{theorem}

Note that the spiral wedge described in (iii) is the image of the
spirallike function with respect to a boundary point
$h_{\nu,a}(z):=\left(\frac{1-z}{1+e^{-2ia}z}\right)^\nu$, so $f$
is subordinate to $h_{\nu,a}$ (see Definition~1.1 below).

In the present work we establish a number of distortion results
and sharp covering theorems for some subclasses of univalent
functions spirallike with respect to a boundary point. Namely, we
consider the following classes of functions.

\begin{defin}\label{defgmb}
We say that a univalent function $f\in\Hol(\Delta,\C),\ f(0)=1$,
belongs to the class $\gmb$ (is $\mu$-spirallike of order $\beta$
with respect to a boundary point), where
$\mu\in\Omega:=\left\{\lambda\in\C:\ |\lambda-1|\le0,\
\lambda\not=0 \right\}$ and $\beta\in[0,1)$, if $f$ satisfies the
inequality:
\be\label{gmb}
\Re\left(\frac2\mu\cdot\frac{zf'(z)}{f(z)}+\frac{1+z}{1-z}\right)>\beta,\quad
z\in\Delta.
\ee
\end{defin}

\section{Integral representation for the class $\gmb$}

\setcounter{equation}{0}

In this section we prove an integral representation for functions
of the class $\gmb$ and deduce some its consequences.

\begin{theorem}\label{P1.1} 
$f\in\gmb$ if and only if $f$ admits the following representation:
\be\label{int}
f(z)=(1-z)^\mu\exp\left[-\mu(1-\beta)\oint\ln(1-z\bar\zeta)d\sigma(\zeta)\right],
\ee
where $\sigma$ is a probability measure on the unit circle
$|\zeta|=1$.
\end{theorem}

\pr First, suppose that the function $f$ is represented by
(\ref{int}). Differentiating~(\ref{int}) and doing a simple
calculation we check that inequality (\ref{gmb}) holds.

Suppose now that $f\in\gmb$. Consider the function
\[
h(z)=\dst\frac1{1-\beta}
\left[\frac2\mu\frac{zf'(z)}{f(z)}+\frac{1+z}{1-z}-\beta\right],
\]
which obviously belongs to the Carath\'eodory class and thus can
be represented by the Riesz--Herglotz formula:
\[
h(z)=\oint\frac{1+z\bar\zeta}{1-z\bar\zeta}\,d\sigma(\zeta),
\]
where $\sigma$ is a probability measure. This implies that
\[
\frac{f'(z)}{\mu f(z)}=\oint_{|\zeta|=1}\left[\frac{(1-\beta)
\bar\zeta}{1-z\bar\zeta} -\frac 1{1-z} \right]d\sigma(\zeta).
\]
Integrating both sides of this equations we obtain~(\ref{int}).
\epr

{\bf Remark.} It follows by Theorem~\ref{AES}, Theorem~\ref{P1.1}
and Definition~\ref{defgmb} that any function $f$ spirallike with
respect to a boundary point has the integral representation
\[
f(z)=(1-z)^\mu\exp\left[-\mu\oint\ln(1-z\bar\zeta)\,
d\sigma(\zeta)\right]
\]
for some $\mu\in\Omega$ and a probability measure $d\sigma$.
Define another measure $d\widetilde{\sigma}$ such that $d\sigma=
\beta d\lambda+(1-\beta)d\widetilde{\sigma}$, where $d\lambda$ is
the normalized Lebesgue measure on the unit circle $\pl\Delta$.
Obviously, $\int\limits_{\pl\Delta}d\widetilde{\sigma}=1$. Since
the integral of the antiholomorphic function $\log(1-z\bar\zeta)$
with respect to the Lebesgue measure $d\lambda$ is zero, we have
that $f$ belongs to the class $\gmb$ if and only if the measure
$d\widetilde{\sigma}=\frac1{1-\beta}\left(d\sigma-
\beta d\lambda\right)$ is positive.

\vspace{2mm}

The following assertion is an immediate consequence of
Theorem~\ref{P1.1}.

\begin{corol}\label{cor1}
Let $\mu_1,\mu_2\in\Omega$ and $\beta_1,\beta_2\in[0,1)$. Let
$f\in\Hol(\Delta,\C)$ be univalent. Then $f\in G(\mu_1,\beta_1)$
if and only if the function
\[
\widetilde{f}(z)=(1-z)^{\mu_2(\beta_2-\beta_1)}
(f(z))^{\frac{\mu_2(1-\beta_2)}{\mu_1(1-\beta_1)}}
\]
belongs to $G(\mu_2,\beta_2)$. In particular, $f\in\gmb$ if and
only if $\widetilde{f}(z)=(f(z))^{\frac1{\mu}}$ belongs to
$G(1,\beta)$.
\end{corol}

\begin{corol}
Let $\mu_1,\mu_2\in\Omega,\,\mu_1=r\mu_2$ with $r\le1$, and let
$\beta_2\le r\beta_1$. Then $G(\mu_1,\beta_1)\subset
G(\mu_2,\beta_2)$.
\end{corol}

\pr Let $f\in G(\mu_1,\beta_1)$. We have
\bep
f(z)&&=(1-z)^{\mu_1}\exp\left[-\mu_1(1-\beta_1)\oint_{|\zeta|=1}
\ln(1-z\bar\zeta)d\sigma(\zeta)\right]\\
&&=(1-z)^{\mu_2}\exp\left[-\mu_2(1-r\beta_1)\oint_{|\zeta|=1}
\ln(1-z\bar\zeta)d\widetilde{\sigma}(\zeta)\right],
\eep
where $d\widetilde\sigma$ is a probability measure on the unit
circle defined by
\[
d\widetilde{\sigma}(\zeta)=\frac{r(1-\beta_1)}{1-r\beta_1}\,d\sigma(\zeta)+
\frac{1-r}{1-r\beta_1}\,\delta(\zeta),
\]
with $\delta$ the Dirac measure at the point $\zeta=1$. By
Theorem~1.1, $f$ satisfies the inequality
\[
\Re\left(\frac2{\mu_2}\cdot\frac{zf'(z)}{f(z)}+\frac{1+z}{1-z}\right)>r\beta_1\ge\beta_2.
\]
Therefore, $f\in G(\mu_2,\beta_2)$. \epr

\begin{corol}\label{dense}
The set
\[
\bigcup^\infty_{n=1}\left\{\left(\frac{1-z}{\prod^n_{j=1}(1-z\bar\zeta_j)^{\lambda_j}}\right)^\mu:\,
\lambda_j\ge0,\ \sum^n_{j=1}\lambda_j =1-\beta, \ |\zeta_j|=1
\right\}
\]
is dense in $\gmb$ in the topology of uniform convergence on
compact subsets of $\Delta$.
\end{corol}

\pr Replacing the integral in (\ref{int}) by approximating sums,
we have for any $f\in G(\mu,\beta)$:
\bep
f(z)&&=\lim_{n\to\infty}(1-z)^\mu\exp\left[-\mu(1-\beta)
\sum_{j=1}^{n} \ln(1-z\bar\zeta_j)\sigma_j\right]\\
&&=\lim_{n\to\infty}(1-z)^\mu\prod_{j=1}^n
(1-z\bar\zeta_j)^{-\mu(1-\beta)\sigma_j},
\eep
with $\sum_{j=1}^{n}\sigma_j=1$. Putting
$\lambda_j=(1-\beta)\sigma_j$, we complete the proof. \epr

To formulate our next result we need the notion of subordination.

\begin{defin}
A function $s_1\in\Hol(\Delta,\C)$ is said to be subordinate to a
function $s_2\in\Hol(\Delta,\C)$ ($s_1\prec s_2$) if there exists
a holomorphic function $\omega$ on $\Delta$ with
$|\omega(z)|\le|z|,\ z\in\Delta$, such that $s_1=s_2\circ\omega$.
\end{defin}

\begin{corol}\label{subord}
Let $f\in\gmb$. Then
$\dst\frac{f(z)}{(1-z)^\mu}\prec\frac1{(1-z)^{\mu(1-\beta)}}\,$.
\end{corol}

\pr Since the function $\log(1-z)$ is convex, for any probability
measure $d\sigma$ there exists an analytic function
$\omega:\Delta\mapsto\Delta$ with $\omega(0)=0$ such that
\[
\oint\ln(1-z\bar\zeta)d\sigma(\zeta)=\log(1-\omega(z)).
\]
Using Theorem~\ref{P1.1} we obtain that
\bep
\frac{f(z)}{(1-z)^\mu}=
\exp\left[-\mu(1-\beta)\oint\ln(1-z\bar\zeta)d\sigma(\zeta)\right]
=\frac1{(1-\omega(z))^{\mu(1-\beta)}}\,,
\eep
which proves our assertion. \epr

It turns out that for the classes $\gmb$ one can improve the
result of Theorem~\ref{ES} (iii).

\begin{corol}
Let $f\in\gmb$. Then $f\prec
h_{\nu,a}(z):=\left(\frac{1-z}{1+e^{-2ia}z}\right)^\nu$, where
$\frac\nu\mu$ is a real number satisfying
$\beta\le\frac\nu\mu\le1$ and $|a|<\frac\pi2\frac\mu\nu(1-\beta)$.
Consequently, $|a|< \frac\pi2\cdot\min\left(1,\frac{1-\beta}\beta
\right)$.
\end{corol}

\pr Using representation (\ref{int}) one can estimate the numbers
$\nu$ and $a$ from Theorem~\ref{ES} (i). Specifically,
\[
\frac{f'(r)(r-1)}{\mu
f(r)}=1-(1-\beta)(1-r)\int_{\pl\Delta}\frac{\bar\zeta}{1-r\bar\zeta}d\sigma(\zeta).
\]
Thus $\frac\nu\mu=\lim\limits_{r\to1^-}\frac{f'(r)(r-1)}{\mu
f(r)}\ge\beta$.

Further,
\[
|a|=\lim\limits_{r\to1^-}=\frac\mu\nu(1-\beta)\left|\int_{\pl\Delta}\arg(1-r\bar\zeta)d\sigma(\zeta)
\right|.
\]
Thus the assertion follows. \epr

\noindent{\bf Remark. }For the case when the number $\nu$ in
(\ref{mu}) (consequently, $\mu$ in (\ref{gmb})) is real,
Theorem~\ref{ES} asserts that the image of a function $f$ is
contained in the wedge of angle $\nu\pi$ with the midline $\arg
w=\nu a$. The last Corollary implies that if $f\in\gmb$ then the
angle can not be less than $\pi\mu\beta$ and the argument of the
midline can not be greater than $\frac\pi2\mu(1-\beta)$.

\section{Estimates for $f$ and $f'$}

\setcounter{equation}{0}

In this section we establish distortion theorems for classes
$\gmb$. In particular, given a point $z\in\Delta$, we find the set
of values for some functionals on this classes. For the class of
Robertson (the class of starlike functions with respect to a
boundary point having image in a half-plane) similar results can
be found in \cite{T} (see also \cite{C-O, SH-SEM, D-L}).

\begin{theorem}\label{P2.1}
For each fixed $z\in\Delta$ we have
\be\label{dist1}
(i)\hspace{1.3cm}\left\{\dst\frac{1-z}{f(z)^{1/\mu}}:\
f\in\gmb\right\}=\left\{ (1+\lambda z)^{1-\beta},\
|\lambda|\le1\right\}\hspace{1.3cm}
\ee
and
\be\label{int_3}
(ii)\hspace{1.3cm}&&\left\{\frac{f'(z)}{\mu
f(z)}+\frac1{1-z}:\,f\in\gmb\right\}\\&&\hspace{2cm}=\left\{w:\,
\left|w-\frac{(1-\beta)\bar z}
{1-|z|^2}\right|\le\frac{1-\beta}{1-|z|^2}
\right\}.\hspace{1.6cm}\nonumber
\ee
Furthermore, if $f\in\gmb$ and $z\in\Delta,\ z\not=0$, one of the
relations
\[
\frac{1-z}{f(z)^{1/\mu}}\in \left\{ (1+\lambda z)^{1-\beta},\
|\lambda|=1\right\}
\]
and
\[
\frac{f'(z)}{\mu f(z)}+\frac1{1-z} \in \left\{w:\,
\left|w-\frac{(1-\beta)\bar z}
{1-|z|^2}\right|=\frac{1-\beta}{1-|z|^2}\right\}
\]
holds only if
\be\label{func}
f(z)=\frac{(1-z)^\mu}{(1-z\bar\xi)^{(1-\beta)\mu}}\,,\quad
\xi\in\pl\Delta.
\ee
\end{theorem}

\pr (i) By Theorem~\ref{P1.1}
\be\label{int_log}
\log\left[\frac{1-z}{f(z)^{1/\mu}}\right]=
(1-\beta)\oint_{\pl\Delta}\log(1-z\bar\zeta)d\sigma(\zeta).
\ee
By the Carath\'eodory principle
\[
\left\{\log\left(\dst\frac{1-z}{f(z)^{1/\mu}}\right):\
f\in\gmb\right\}={\rm Conv} \left\{(1-\beta)\log(1-z\bar\zeta),\
\zeta\in\pl\Delta \right\},
\]
where ${\rm Conv}$ denotes the convex hull.

Since for each $z\in\Delta$ the function
$g(w):=(1-\beta)\log(1-zw)$ maps $\Delta$ onto a strictly convex
domain, this formula and (\ref{int_log}) imply that the value
$\log\left[\frac{1-z}{f(z)^{1/\mu}}\right]$ belongs to the image
$g(\Delta)$ except for the case when the measure $d\sigma(\zeta)$
is the Dirac $\delta$-function at some boundary point $\xi$. Hence
the assertion follows.

(ii) Once again using representation (\ref{int}) we get
\[
\frac{f'(z)}{\mu
f(z)}+\frac1{1-z}=(1-\beta)\oint_{\pl\Delta}\frac{\bar\zeta}{1-z\bar\zeta}\,d\sigma(\zeta).
\]
For each fixed $z\in\Delta$ the function $g:\,\pl\Delta\mapsto\C$
defined by $g(\zeta):=\frac{\bar\zeta}{1-z\bar\zeta}$ maps the
unit circle onto the circle $\dst\left\{w:\ \left|w-\frac{\bar
z}{1-|z|^2}\right|= \frac1{1-|z|^2} \right\}\,$. Therefore,
(\ref{int_3}) follows from the Carath\'eodory principle. \epr

The following two assertions are immediate consequences of
Theorem~\ref{P2.1}~(i).
\begin{corol}
$\bigcup\limits_{z\in\Delta}\left\{\dst\frac{1-z}{f(z)^{1/\mu}}:\
f\in\gmb\right\}=\left\{ (1+\lambda)^{1-\beta},\
|\lambda|\le1\right\}$.
\end{corol}

\begin{corol}
For each $z\in\Delta$ and $f\in\gmb$,
\[
(1-|z|)^{1-\beta}\le\left|\frac{1-z}{f(z)^{1/\mu}}\right|\le(1+|z|)^{1-\beta}
\]
and
\[
\left|\arg\left(\frac{1-z}{f(z)^{1/\mu}}\right)\right|\le(1-\beta)\arcsin|z|,
\]
where for $z\not=0$ equality is achieved only for the functions
(\ref{func}) at the points $z=\pm|z|\xi$ and $z=|z|\xi e^{\pm
i\arccos|z|}$, respectively. In particular, if $\mu$ is real, then
\[
\frac{|1-z|^\mu}{(1+|z|)^{\mu(1-\beta)}}\le\left|f(z)\right|\le
\frac{|1-z|^\mu}{(1-|z|)^{\mu(1-\beta)}}
\]
\end{corol}

In the case when $\mu\in(0,2]$ (i.e., $f\in\gmb$ is, in fact,
starlike) one obtains the following estimate for $f'$.

\begin{corol}
For each $z\in\Delta$ and $f\in\gmb,\ \mu\in(0,2],\
\beta\in[0,1)$,
\bep
\frac{\mu|1-z|^\mu}{(1-|z|^2)(1+|z|)^{\mu(1-\beta)}}
\left(\left|\frac{1-\bar z}{1-z}+\beta\bar z
\right|-1+\beta\right)\le|f'(z)|\le\\  \le
\frac{\mu|1-z|^\mu}{(1-|z|^2)(1-|z|)^{\mu(1-\beta)}}
\left(\left|\frac{1-\bar z}{1-z}+\beta\bar z \right|+1-\beta
\right).
\eep
In particular,
\[
|f'(z)|\le\frac{2\mu|1-z|^\mu}{(1-|z|^2)(1-|z|)^{\mu(1-\beta)}}\,.
\]
\end{corol}
The proof follows from Theorem~\ref{P2.1} (ii) and Corollary~2.2.

\section{Covering theorem}

\setcounter{equation}{0}

In this section we prove a covering theorem for classes $\gmb$ of
spirallike functions with respect to a boundary point. Note that
if for some class of functions there exists a domain $D$ such that
$D\subset f(\Delta)$ for each function $f$ of the class, and the
function $f_0\in\Hol(\Delta,\C)$ maps the open unit disk onto $D$
conformly, then the corresponding covering result can be expressed
as the subordination $f_0\prec f$.

\begin{theorem}
Let $f\in\gmb$. Denote $f_0(z)=(1-z)^{\mu\beta}$. Then $f_0\prec
f$. This subordination is sharp since $f_0\in\gmb$.
\end{theorem}

\pr By Corollary 1.1 it is sufficient to prove this assertion for
the case $\mu=1$. First we suppose that $f\in\gub$ has the form
\[
f(z)=\frac{1-z}{\prod^n_{j=1}(1-z\overline{\zeta_j})^{(1-\beta)\sigma_j}},\
\sum^n_{j=1}\lambda_j =1-\beta, \ |\zeta_j|=1.
\]

Let $z_0=e^{i\phi}\in\pl\Delta,\
z_0\not=\zeta_1,\zeta_2,\ldots,\zeta_n\,$. Consider the two
univalent convex functions
\[
g_0(z)=\log f_0(z)=\beta\log(1-z)
\]
and
\[
g(z)=\log(1-z_0)-(1-\beta)\log(1-z_0z).
\]
It is easy to see that $g_0(z_0)=g(1)$, i.e., this is their common
boundary point. Furthermore,
\[
\left.zg_0'(z)\right|_{z=z_0}=\frac{-\beta
z_0}{1-z_0}\quad\mbox{and}\quad
\left.zg'(z)\right|_{z=1}=\frac{(1-\beta)z_0}{1-z_0}\,,
\]
i.e., the exterior normal vectors to the images $g_0(\Delta)$ and
$g(\Delta)$ at this common point have reverse directions. Because
these functions are convex,\linebreak ${g_0(\Delta)\bigcap
g(\Delta)=\emptyset}$. Thus,
\[
\log f(z_0)=\sum^n_{j=1} \sigma_j\left[\log(1-z_0)-
(1-\beta)\log(1-z_0\overline{\zeta_j})\right]
\]
belongs to $\overline{g(\Delta)}$ and does not belong to
$g_0(\Delta)=\log f_0(\Delta)$. Since the point $f(z_0)$ is an
arbitrary finite boundary point of the image $f(\Delta)$, we have
proved that the boundary $\pl f(\Delta)$ does not intersect
$f_0(\Delta)$. This fact implies that $f_0$ is subordinate to any
function $f\in\gub$ that has the form described above.

In light of Corollary~\ref{dense} and by the Carath\'eodory Kernel
Convergence Theorem we deduce the assertion of the theorem. \epr

\noindent{\bf Example.} Consider the function
\[
f(z)=\frac{1-z}{(1-0.9z-0.4iz)^{0.2}(1-0.9z+0.4iz)^{0.2}}\,.
\]
A simple calculation shows that $f\in\gmb$ with $\mu=1$ and
$\beta=0.6$. Hence by Theorem~3.1, its image $f(\Delta)$ contains
the image $f_0(\Delta)$, where $f_0(z)=(1-z)^{0.6}$ (see Figure).

\begin{center}
\begin{figure}[h]
\centering
\includegraphics[angle=270,width=6cm,totalheight=6cm]{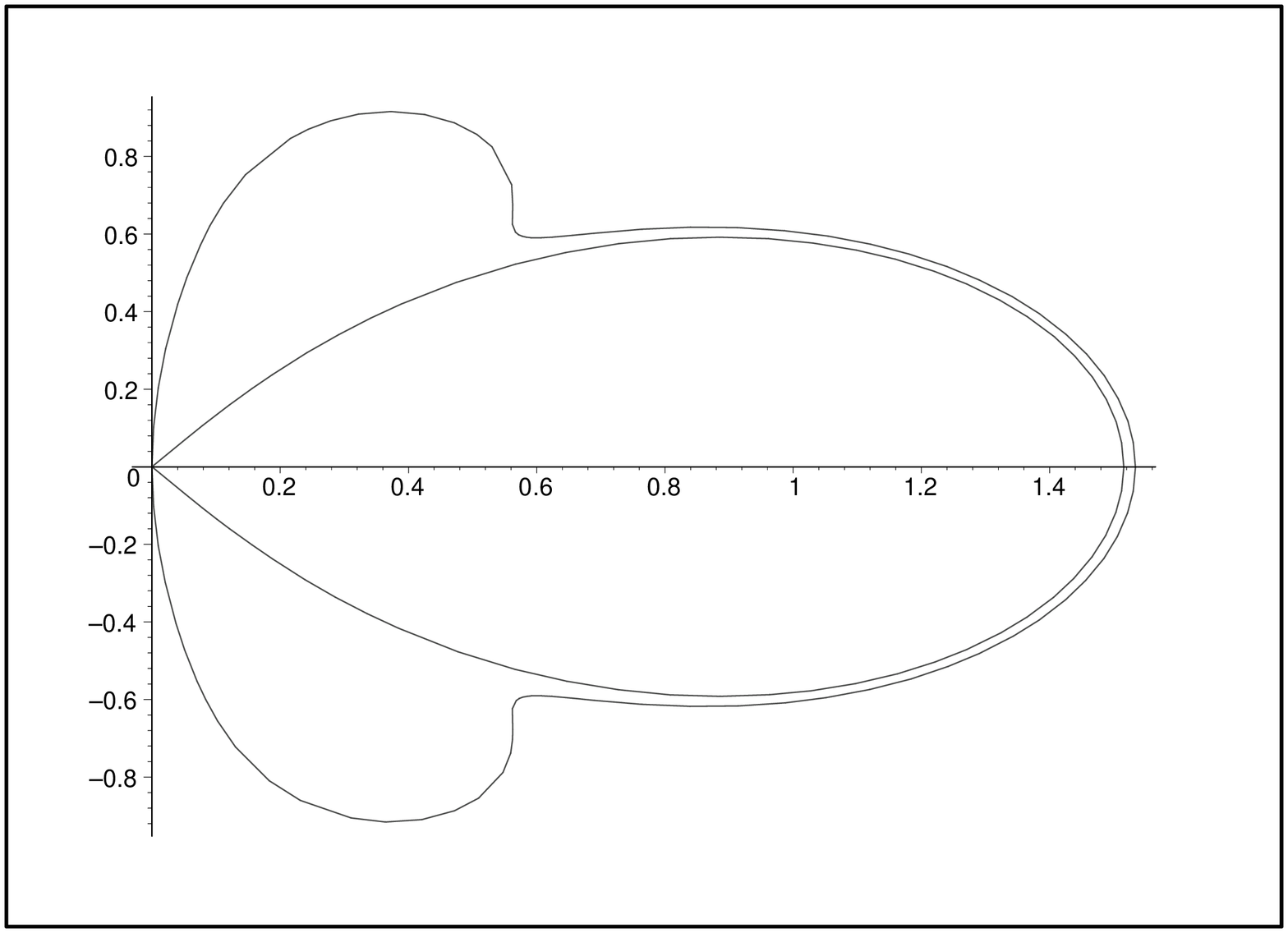}
\end{figure}
\end{center}

\begin{corol}
If $f\in\gmb$ with $\mu\in(0,2]$ and $\beta\in(0,1)$ then the
image of $f$ covers the disk $\left\{w:\, |w-1|<r(\mu\beta)
\right\} $, where
\[
r(s)=\left\{ \begin{array}{c}
\sqrt{1+2^{2s}-2^{s+1}},\quad\mbox{when }s\in(0,1]
\\ 1,\hspace{2.6cm}   \mbox{when }s\in(1,2]
\end{array}\right. .
\]
\end{corol}

\pr It suffices to find the distance between the point $z=1$ and
the curve\linebreak $\{(1-z)^s,\ z\in\pl\Delta \}$, where
$s=\mu\beta$. Indeed, setting $z=-e^{it},\ t\in(-\pi,\pi)$,
consider the function
\bep
a_s(t)=\left|(1+e^{it})^s-1 \right|^2=\left(2\cos\frac
t2\right)^{2s}+1-2\left(2\cos\frac t2\right)^{s}\cos\frac{st}2
\eep
This function is even in $t$, $a_s$ increases on the segment $0\le
t\le\pi$ when $s<1$, and $a_s$ decreases on the same segment when
$s>1$. So,
\[
\min_{-\pi<t<\pi}a_s(t)=\left\{ \begin{array}{c}
a_s(0),\quad\mbox{when }s\in(0,1]
\\ a_s(\pi),\quad\mbox{when }s\in(1,2]
\end{array}\right. .
\]
Calculating these values we complete the proof. \epr

\noindent{\bf Remark.} In \cite{C-O} Chen and Owa have proved a
covering theorem for starlike functions with respect to a boundary
point. Using our notation their result can be quoted as follows:
{\it If $f\in\gmb$ for some $\mu,\beta\in[0,1]$ then its image
covers the disk $\left\{w:|w-1|<\frac{\mu\beta}4\right\}$.} The
radius of the covered disk which was found in the last Corollary
is considerably larger than one due to Chen and Owa. For example,
if $f\in\gmb$ with $\mu\beta=1$, then by the Theorem of Chen and
Owa $f(\Delta)$ covers the disk of radius $1/4$; our result
asserts that the disk covered is of radius $1$.

\vspace{2mm}

The following assertion connects the classes $\gmb$ with
spirallike functions with respect to an interior point.

\begin{propo}\label{P1.2}
Let $\beta\in[0,1)$, and let a complex number $\mu=r\,e^{i\phi}$
belong to $\Omega=\left\{\lambda:\ |\lambda-1|\le1,\
\lambda\not=0\right\}$. Then a univalent function $f,\,f(0)=1,$
belongs to the class $\gmb$ if and only if the function $s$
defined by
\be\label{s}
s(z):=\frac{zf(z)}{(1-z)^\mu}
\ee
is $\phi$-spirallike of order $\dst
\cos\phi-\frac{r(1-\beta)}2\,$, i.e., $s$ satisfies the condition
\be\label{s1}
\Re\left(e^{-i\phi}\frac{zs'(z)}{s(z)}\right)>\cos\phi-\frac{r(1-\beta)}
2\,,\quad z\in\Delta.
\ee
\end{propo}

\noindent{\bf Remark.} Note that the relation
$\mu=r\,e^{i\phi}\in\Omega$ implies that \linebreak
$\cos\phi-\frac{r(1-\beta)}2\ge0$, and then by a result of
\v{S}pa\v{c}ek (see \cite{GAW}) the function $s$ is univalent in
$\Delta$.

\vspace{2mm}

\pr Let a function $s$ be defined by (\ref{s}). Then
\bep
&& e^{-i\phi}\frac{zs'(z)}{s(z)}=\frac r\mu
\left(1+\frac{zf'(z)}{f(z)}+\frac{\mu z}{1-z}\right)\\ && =\frac
r2 \left(\frac{2zf'(z)}{\mu
f(z)}+\frac{1+z}{1-z}-\beta\right)+e^{-i\phi}-\frac{r(1-\beta)}2\,.
\eep
This equality implies the assertion. \epr

Exactly, as in \cite{A-E-S}, using Proposition 3.1 and results of
Ruscheweyh \cite{RuSt-76}, one can conclude the following.

\begin{corol}
Let $f:\Delta\mapsto\C$ be a holomorphic function and $f(0)=1$.
Let $\mu=re^{i\phi}\in\Omega$ and $\beta\in[0,1)$. Then $f\in\gmb$
if and only if one of the following conditions holds:

(a) for all $u,v\in\overline\DD$
\[
\displaystyle\left(\frac{1-uz}{1-vz}\right)^\mu\,\frac{f(vz)}{f(uz)}\prec
\left(\frac{1-uz}{1-vz}\right)^{(1-\beta)\mu};
\]

(b) for all $t\in(0,2\cos\phi)$
\[
\displaystyle\left|\frac{f(z(1-e^{-i\phi}t))}{f(z)}\right|\leq
\left|\left(\frac{1-z(1-e^{-i\phi}t)}{1-z}\right)^\mu\right|\cdot
\left(1-\frac{t}{2\cos\phi}\right)^{-\Re\mu(1-\beta)}.
\]
\end{corol}
(We omit the proof since it repeats one in \cite{A-E-S}.)

Setting $u=0$ and $v=1$ in Corollary 3.2, we obtain an alternative
proof of Corollary 1.4 above.

The next assertion follows from Theorem 3.1 and
Proposition~\ref{P1.2}.
\begin{corol}
Let $\phi\in\left(-\frac\pi2,\frac\pi2\right),\ \alpha<\cos\phi,$
and $\beta\in\left(0,\frac\alpha{\cos\phi}\right]$. Denote
${\mu=e^{i\phi}\frac{2(\cos\phi-\alpha)}{1-\beta}}$. For each
univalent function $s$ normalized by $s(0)=0$, ${s'(0)=1}$ and
satisfying
${\dst\Re\left(e^{-i\phi}\frac{zs'(z)}{s(z)}\right)>\alpha,}$ the
image of the function
\[
1-(1-z)^\frac1\beta\left(\frac{s(z)}z\right)^\frac1{\mu\beta}
\]
covers the open unit disk, and, consequently, the image of the
function
\[
\frac{(1-z)^\frac1\beta\left(\frac{s(z)}z
\right)^\frac1{\mu\beta}}{2-(1-z)^\frac1\beta\left(\frac{s(z)}z
\right)^\frac1{\mu\beta}}
\]
covers the right half-plane.
\end{corol}

In particular, setting $\phi=0,\ \alpha=\beta=\frac12$, we get
that for each function $s$ starlike of order $\frac12$ the image
of the function $1-(1-z)^2\left(\frac{s(z)}z\right)$ contains the
open unit disk, and the image of the function $\dst\frac{(1-z)^2
s(z)}{2z-(1-z)^2 s(z)}$ contains the right half-plane.

\vspace{2mm}

\noindent\textbf{ACKNOWLEDGMENT.} The author is grateful to David
Shoikhet for useful discussions and comments.


\vspace{2mm}

\begin{thebibliography}{999}

\bibitem{A-A-S} A. S. Abdullah, R. M. Ali and V. Singh,
On functions starlike with respect to a boundary point, {\it Ann.
Univ. Mariae Curie-Sk\l odowska Sect. A} {\bf 50} (1996), 7--15.

\bibitem{A-E-S} D. Aharonov, M.~Elin and D.~Shoikhet, Spirallike functions with
respect to a boundary point, {\it J. Math. Anal. Appl.}  {\bf 280}
(2003), 17--29.

\bibitem{C-O} M. Chen and S. Owa,
Generalization of Robertson functions, in: {\it Complex analysis
and its applications (Hong Kong, 1993)}, 159--165, Pitman Res.
Notes Math. Ser., 305, Longman Sci. Tech., Harlow, 1994.

\bibitem{D-L} J. Dong and L. Liu,
Distortion properties of a class of functions starlike with
respect of a boundary point, {\it Heilongjiang Daxue Ziran Kexue
Xuebao} {\bf 15} (1998), 1--6.

\bibitem{E-R-S1} M. Elin, S. Reich and D. Shoikhet,
Dynamics of inequalities in geometric function theory,
{\it Journ. of Inequal. $\&$ Appl.} {\bf 6} (2001), 651--664.

\bibitem{E-R-S2} M. Elin, S. Reich and D. Shoikhet,
Holomorphically accretive mappings and spiral-shaped functions of
proper contractions, {\it Nonlinear Analysis Forum} {\bf 5}
(2000), 149--161.

\bibitem{E-S2} M. Elin and D. Shoikhet, Univalent functions of proper
contractions spirallike with respect to a boundary point, in: {\it
Multidimensional Complex Analysis} Krasnoyarsk, (2002), 28--36.

\bibitem{GAW} A. W. Goodman,{\it Univalent Functions}, Vols. I, II,
Mariner Publ. Co., Tampa, FL, 1983.

\bibitem{LeA} A. Lecko,
On the class of functions starlike with respect to a boundary
point, {\it J. Math. Anal. Appl.} {\bf 261} (2001), 649--664.

\bibitem{LeA1} A. Lecko,
The class of functions spirallike with respect to a boundary
point, {\it Int. J. Math. Math. Sci.} (2004) 2133--2143.

\bibitem{LA-LA} A. Lecko and A. Lyzzaik,
A note on univalent functions starlike with respect to a boundary
point, {\it J. Math. Anal. Appl.} {\bf 282} (2003), 846--851.

\bibitem{LA} A. Lyzzaik, On a conjecture of M.~S.~Robertson,
{\it Proc. Amer. Math. Soc.} {\bf 91} (1984), 108--110.

\bibitem{RMS} M. S. Robertson, Univalent functions starlike with respect
to a boundary point, {\it J. Math. Anal. Appl.} {\bf 81} (1981),
327--345.

\bibitem{RuSt-76} S. Ruscheweyh, A subordination theorem for $\Phi$-like
functions, {\it J. London Math. Soc.} {\bf 13} {\it 2} (1976),
275--280.

\bibitem{T} P. Todorov,
On the univalent functions starlike with respect to a boundary
point, {\it Proc. Amer. Math. Soc.} {\bf 97} (1986), 602--604.

\bibitem{SD} D. Shoikhet, {\it Semigroups in Geometrical Function Theory},
Kluwer Acad. Publ., 2001.

\bibitem{SD1} D. Shoikhet,
Representations of holomorphic generators and distortion theorems
for spirallike functions with respect to a boundary point, {\it
Int. J. Pure Appl. Math.} {\bf 5} (2003), 335--361.

\bibitem{SH-SEM} H. Silverman and E. M. Silvia, Subclasses of univalent
functions starlike with respect to a boundary point, {\it Houston
J. Math.} {\bf 16} (1990), 289--299.

\bibitem{Z-Q} Y. L. Zhang and F. C. Qian,
Extremal problems for the class of Robertson functions, {\it Acta
Math. Sinica} {\bf 33} (1990), 601--609.

\end{thebibliography}
\end{document}